\documentclass[11pt]{article}   
\usepackage{amsthm, amsmath}   
\usepackage{latexsym,amssymb}
\usepackage{algorithmic}
\usepackage{pst-node}
\usepackage{amssymb}

\def\QED{\hfill$\Box$}

\def\P{{\mathbb P}}
\def\N{{\mathbb N}}
\def\A{{\mathbb A}}
\def\Z{{\mathbb Z}}
\def\Q{{\mathbb Q}}

\def\gcd{{\rm{gcd}}}
\def\lcm{{\rm{lcm}}}
\def\mA{\mathcal A}
\def\mB{\mathcal B}
\def\mC{\mathcal C}
\def\mS{\mathcal S}
\def\omA{\mA^{\star}}

\def\omC{\mC^{\star}}
\def\IA{I_{\mathcal A}}
\def\IB{I_{\mathcal B}}
\def\oIA{I_{\mA^{\star}}}
\def\oIB{I_{\mB^{\star}}}

\def\gcd{{\rm{gcd}}}
\def\lcm{{\rm{lcm}}}
\def\g{{\rm g}}

\def\x{{\mathbf{x}}}
\def\kx{k[\x]}
\def\t{{\mathbf{t}}}
\def\kt{k[\t]}

\theoremstyle{plain}

\newtheorem{Theorem}{Theorem}[section]
\newtheorem{Lemma}[Theorem]{Lemma}
\newtheorem{Corollary}[Theorem]{Corollary}
\newtheorem{Proposition}[Theorem]{Proposition}
\newtheorem{Remark}[Theorem]{Remark}

\newtheorem{theorem}{Theorem}

\theoremstyle{definition}

\newtheorem{example}[theorem]{Example}

\title{Complete intersections in certain affine and projective monomial curves \footnote{Partially supported by Ministerio
de Ciencia e Innovaci\'on - Espa\~na (MTM2010-20279-C02-02)}}

\author{Isabel Bermejo, Ignacio Garc\'{\i}a-Marco}

\begin{document}
\maketitle

\begin{abstract}Let $k$ be an arbitrary field, the purpose of this work
is to provide families of positive integers $\mA =
\{d_1,\ldots,d_n\}$ such that either the toric ideal $\IA$ of the
affine monomial curve $\mathcal C = \{(t^{d_1},\ldots,\,t^{d_n}) \ |
\ t \in k\} \subset \A_k^n$ or the toric ideal $\oIA$ of its
projective closure $\omC \subset \P_k^n$ is a complete intersection.
More precisely, we characterize the complete intersection property
for $\IA$ and for $\oIA$ when:
\begin{itemize}\item[(a)] $\mA$ is a generalized arithmetic
sequence, \item[(b)] $\mA \setminus \{d_n\}$ is a generalized
arithmetic sequence and $d_n \in \Z^+$,
\item[(c)] $\mA$ consists of certain terms of the $(p,q)$-Fibonacci
sequence, and \item[(d)] $\mA$ consists of certain terms of the
$(p,q)$-Lucas sequence.\end{itemize} The results in this paper arise
as consequences of those in \cite{BGS, BGsimplicial} and some new
results regarding the toric ideal of the curve.
\end{abstract}

\section{Introduction}
Let $k$ be an arbitrary field and $\kx= k[x_1,\ldots,x_n]$ and $\kt
= k[t_1,\ldots,t_m]$ two polynomial rings over $k$. A {\it binomial}
$f$ in $\kx$ is a difference of two monomials. Let $\mB
=\{b_1,\ldots,b_n\}$ be a set of nonzero vectors in $\N^m$; each
vector $b_i = (b_{i1},\ldots,b_{im})$ corresponds to a monomial
$\t^{b_i} = t_1^{b_{i1}} \cdots t_m^{b_{im}}$ in $\kt$.

\medskip The {\it toric ideal} determined by $\mB$ is the  kernel of
the homomorphism of $k$-algebras $$\varphi\colon \kx \rightarrow
\kt;\ x_i\longmapsto \t^{b_i}$$ and is denoted by
 $\IB$. By \cite[Corollary~4.3]{Sturm}, $\IB$ is an $\mB$-{\it homogeneous binomial\/} ideal, i.e.,
if one sets the $\mB$-degree of a monomial $\x^{\alpha}  \in \kx$ as
${\rm deg}_{\mB}(\x^{\alpha}) := \alpha_1 b_1 + \cdots + \alpha_n
b_n \in \N^m$, and says that a polynomial $f \in k[x]$ is
$\mB$-homogeneous if its monomials have the same $\mB$-degree, then
$\IB$ is generated by $\mB$-homogeneous binomials. According to
\cite[Lemma 4.2]{Sturm}, the height of $\IB$ is ${\rm ht}(\IB) =
n-{\rm rk}(\Z \mB)$, where ${\rm rk}(\Z \mB)$ denotes the rank of
the subgroup of $\Z^m$ generated by $\mB$. The ideal $\IB$ is a {\it
complete intersection} if $\mu(\IB) = {\rm ht}(\IB)$, where
$\mu(\IB)$ denotes the minimal number of generators of $\IB$.
Equivalently, $\IB$ is a complete intersection if there exists a set
of $s = n-{\rm rk}(\Z \mB)$ $\mathcal B$-homogeneous binomials
$g_1,\ldots,g_{s}$ such that $\IA=(g_1,\ldots,g_{s}).$

\smallskip Let $\mA = \{d_1,\ldots,d_n\}$ be a subset of the positive
integers and consider the {\it affine monomial curve} $\mathcal C$
given parametrically by $x_1 = t^{d_1},\ldots,x_n = t^{d_n}$, i.e.,
$$
\mC = \{(t^{d_1},\ldots,t^{d_n})\in\mathbb{A}_k^n\,\vert\, t\in
k\}\,.
$$
By \cite[Corollary~7.1.12]{monalg}, if $k$ is an infinite field the
ideal $I(\mathcal C)$ of polynomials vanishing on $\mathcal C$ is
equal to $\IA$, the toric ideal determined by $\mA$. $\IA$ is called
the {\it toric ideal of} $\mathcal C$.

\smallskip Set $d := {\rm max}(\mA)$ and consider $\omA := \{a_1, \ldots, a_{n-1}, a_n, a_{n+1}\} \subset \N^2,$
where $a_i := (d_i,d - d_i)$ for all $i \in \{1,\ldots,n\}$ and
$a_{n+1} := (0,d)$, and the {\it projective monomial curve $\omC$}
obtained as the projective closure of $\mC$, which is given
parametrically by $x_1 = t^{d_1} u^{d-d_1},\ldots,x_n = t^{d_n} u
^{d-d_n},x_{n + 1} = u^d$, i.e.,
$$
\omC = \{(t^{d_1}u^{d-d_1}:\cdots:t^{d_n}u^{d - d_n}:u^d)\in
\P_k^n\,\vert\, (t:u) \in \P_k^1\}\,.
$$

\noindent Again by \cite[Corollary~7.1.12]{monalg}, if $k$ is an
infinite field, then $I(\omC) = \oIA$. The toric ideal $\oIA$ is
homogeneous, indeed $\oIA$ is the homogenization of $\IA$ with
respect to the variable $x_{n+1}$. $\oIA$ is called the {\it toric
ideal of} $\omC$.

\smallskip
Both $\IA$ and $\oIA$ have height $n-1$; thus $\IA$ (resp. $\oIA$)
is a {\it complete intersection\/} if there exists a system of
$\mA$-homogeneous (resp. homogeneous) binomials $g_1,\ldots,g_{n-1}$
such that $\IA=(g_1,\ldots,g_{n-1})$ (resp. $\oIA =
(g_1,\ldots,g_{n-1}) ).$ Clearly, $\IA$ is a complete intersection
whenever $\oIA$ is, the converse is not true in general.

\smallskip
The aim of this work is to provide families of positive integers
$\mA = \{d_1,\ldots,d_n\}$ such that either the toric ideal $\IA$ or
$\oIA$ is a complete intersection. The starting point are the papers
by Garc\'{\i}a-S\'anchez and Rosales \cite{GR99}, Maloo and Sengupta
\cite{MS} and Fel \cite{Fel}.

\smallskip
In the first one, the authors prove that if $\mA$ is a set of
consecutive positive integers, then $\IA$ is a complete intersection
if and only if $n = 2$ or $n = 3$ and $d_1$ is even. In the second
one, the authors obtain the same characterization when $\mA$ is an
arithmetic sequence provided $\gcd(\mA) = 1$. In Theorem \ref{gen}
we generalize this result to {\it generalized arithmetic sequences}.
We recall that $\mA$ is a generalized arithmetic sequence if there
exists $h \in \Z^+$ such that $\{h d_1,d_2,\ldots,d_n\}$ is an
increasing arithmetic sequence.

\smallskip
Maloo and Sengupta also study the case in which $\mA$ is an
almost-arithmetic sequence, i.e., $\mA \setminus \{d_n\}$ is an
arithmetic sequence and $d_n \in \Z^+$, and prove that $n \leq 4$
provided
 $\IA$ is a complete intersection. In Theorem \ref{almost} we go
 further and characterize when $\IA$ is a complete intersection whenever $\mA
\setminus \{d_n\}$ is a generalized arithmetic sequence and $d_n \in
\Z^+$, where $n \geq 4$.

\smallskip In the third one, the author provides certain conditions so
that the semigroup generated by $\mA = \{d_1,d_2,d_3\}$ where
$\gcd\{d_1,d_2,d_3\} = 1$ and $d_1,d_2,d_3$ are members of the
Fibonacci or Lucas sequence is symmetric. Recall that for a set $\mA
= \{d_1,\ldots,d_n\}$ with $\gcd(\mA) = 1$, setting $\mS := \sum_{i
= 1}^n \N\, d_i$, then the complement $\mS$ in $\N$ is finite and
the largest integer not belonging to $\mS$ is called the {\it
Frobenius number of $\mS$} and denoted by $\g(\mS)$. Moreover, the
semigroup $\mS$ is {\it symmetric} if for every $d \in \Z$, either
$d \in \mS$ or $\g(\mS) - d \in \mS$. It is a classical result due
to Herzog \cite[Theorem 3.10]{Herzog} that whenever $\mA =
\{d_1,d_2,d_3\} \subset \Z^+$ with $\gcd\{d_1,d_2,d_2\} = 1$, then
$\mS$ is symmetric if and only if $\IA$ is a complete intersection.
 In this work we characterize the complete intersection property
for $\IA$ when $\mA$ is a certain subset of either the
$(p,q)$-Fibonacci sequence or the $(p,q)$-Lucas sequence, where $p,q
\in \Z^+$ are relatively prime. We recall that the $(p,q)$-Fibonacci
sequence, denoted by $\{F_n\}_{n\in \N}$, is defined as follows
\begin{center}$F_0 = 0,\,F_1 = 1$ and $F_{n+2} = p\, F_{n+1} + q\, F_n$ for every $n \geq 0$,
\end{center}
and the $(p,q)$-Lucas sequence, denoted by $\{L_n\}_{n \in \N}$, is
defined as \begin{center} $L_0 = 2,\,L_1 = p$ and $L_{n+2} = p\,
L_{n+1} + q\, L_n$ for every $n \geq 0$.
\end{center}
These sequences are natural generalizations of the usual Fibonacci
and Lucas sequences, now called $(1,1)$-Fibonacci sequence and
$(1,1)$-Lucas sequence. In Theorems \ref{fibonacci} and \ref{lucas}
we characterize when $\IA$ is a complete intersection, being $\mA =
\{d_1,\ldots,d_n\}$ with:
\begin{enumerate}
\item[(a)] $d_i = F_{e_i}$, where
$\{e_1,\ldots,e_n\}$ is a generalized arithmetic sequence, and
\item[(b)] $d_i = L_{e_i}$, where
$\{e_1,\ldots,e_n\}$ is an arithmetic sequence.
\end{enumerate}

\smallskip
Moreover, we characterize algorithmically when the toric ideal of a
projective monomial curve is a complete intersection with Algorithm
CI-projective-monomial-curve (see Table
\ref{algoritmocurvamonomialproyectiva}). Using this algorithm we are
able to characterize in Theorems \ref{genproy}, \ref{almostN=4proy},
\ref{fibonacciproy} and \ref{lucasproy}, which are the projective
versions of Theorems \ref{gen}, \ref{almost}, \ref{fibonacci} and
\ref{lucas} respectively, when $\oIA$ is a complete intersection
when $\mA$ belongs to any of the four families already described.

\smallskip
Whenever $\IA$ or $\oIA$ is a complete intersection, we also obtain
a minimal set of generators of the toric ideal. Furthermore, when
$\IA$ is a complete intersection, using the formula described in
\cite[Remark 11]{D} and \cite[Remark 4.5]{BGRV}, the Frobenius
number $\g(\mathcal S)$ of the numerical semigroup $\mathcal S :=
\sum_{i = 1}^n \N \, (d_i / e)$, where $e := \gcd(\mA)$, is also
provided. Indeed, the formula asserts that if $\IA =
(g_1,\ldots,g_{n-1})$ where $g_i$ is $\mA$-homogeneous for all $i
\in \{1,\ldots,n-1\}$, then $ \g(\mathcal S) =
\left(\sum_{i=1}^{n-1} {\rm deg}_{\mA}(g_i) - \sum_{i = 1}^n d_i
\right) / e$.

\smallskip
The results obtained in this work arise as consequences of our
papers \cite{BGS, BGsimplicial} and some new results concerning the
toric ideal of the curves. In \cite{BGS} we exploited the
combinatorial-arithmetical structure of complete intersections given
by the existence of a certain {\it binary tree labeled by}
$\{d_1,\ldots,d_n\}$ stated in \cite[Theorem~4.3]{BGRV} to obtain an
algorithm that determines whether the toric ideal of an affine
monomial curve is a complete intersection. In \cite{BGsimplicial} we
provided some new results concerning complete intersection toric
ideals in general and apply them in order to obtain algorithms that
characterize when either a simplicial toric ideal or a homogeneous
simplicial toric ideal is a complete intersection. Since both $\IA$
and $\oIA$ are simplicial toric ideals and, moreover, $\oIA$ is
homogeneous, the algorithms obtained in \cite{BGsimplicial} apply to
them. These algorithms have been implemented in ANSI C programming
language and also in the distributed library {\tt cisimplicial.lib}
\cite{BGlib2} of {\sc Singular} \cite{DGPS}.

\section{Complete intersection toric ideals associated to affine and projective monomial curves}

This section is devoted to present some new results concerning the
complete intersection property for toric ideals associated to either
affine or projective monomial curves. The results of this section
arise after applying some results of \cite{BGsimplicial} to the
context of affine and projective monomial curves. The main results
of this section are namely Proposition \ref{n-1dist} and Theorem
\ref{coralg}. On one hand, Proposition \ref{n-1dist} provides, under
certain hypothesis, a necessary and sufficient condition for the
toric ideal of an affine monomial curve to be a complete
intersection. This result will be useful in Sections $3, 4$ and $5$.
On the other hand, Theorem \ref{coralg} is a particularization of
\cite[Corollary 5]{BGsimplicial} for toric ideals associated to
projective monomial curves. This result, together with Remark
\ref{notaproyred}, yields Algorithm CI-projective-monomial-curve of
Table \ref{algoritmocurvamonomialproyectiva}, an algorithm for
checking whether the toric ideal of a projective monomial curve is a
complete intersection. This algorithm will be useful in Section $6$.

\medskip
In order to present the new results, we begin our explanation by
briefly describing some results of \cite{BGsimplicial}. For every
set $\mB = \{b_1,\ldots,b_n\}$ of nonzero vectors of $\N^m$, we have
the following results.

\begin{Lemma}\cite[Lemmas 2.1 and 2.2]{BGsimplicial}\label{facilred} \begin{itemize} \item If $b_i \notin \sum_{j \in \{1,\ldots,n\} \atop j \neq i} \Q\, b_j$, then $\IB =
I_{\mB \setminus \{b_i\}} \cdot \kx$. Moreover, $\IB$ is a complete
intersection $\Longleftrightarrow \, I_{\mB \setminus \{b_i\}}$ so
is.
\item If $b_i = \sum_{j \in \{1,\ldots,n\} \atop j \neq i} \alpha_j
b_j \in \sum_{j \in \{1,\ldots,n\} \atop j \neq i} \N\, b_j$, then
$\IB = I_{\mB \setminus \{b_i\}} \cdot \kx + (x_i - \prod_{j \in
\{1,\ldots,n\} \atop j \neq i} x_j^{\alpha_j})$. Moreover, $\IB$ is
a complete intersection $\Longleftrightarrow \, I_{\mB \setminus
\{b_i\}}$ so is.
\end{itemize}
\end{Lemma}

\bigskip Whenever $b_i \in \sum_{j \in \{1,\ldots,n\} \atop j\neq i}
\Q\, b_j$ we denote $B_i := {\rm min}\{B \in \Z^+\, \vert \, B b_i
\in \sum_{j \in \{1,\ldots,n\} \atop j \neq i} \Z b_j\}$ and have
the following result.

\begin{Proposition}\cite[Proposition 2.3]{BGsimplicial}\label{red}
Assume that $b_i \in \sum_{j \in \{1,\ldots,n\} \atop j\neq i} \Q\,
b_j$ for some $i \in \{1,\ldots,n\}$ and set $\mB' :=
\{b_1,\ldots,B_i b_i, \ldots, b_n\}$ and $\rho: \kx \longrightarrow
\kx$ the morphism induced by $\rho(x_i) = x_i^{B_i}$, $\rho(x_j) =
x_j$ for every $j \neq i$. Then, $\IB = \rho(I_{\mB'}) \cdot \kx$.
Moreover, $\IB$ is a complete intersection $\Longleftrightarrow$
$I_{\mB'}$ is a complete intersection.
\end{Proposition}

Applying Lemma \ref{facilred} and Proposition \ref{red} iteratively,
we can associate to $\mB$ a unique subset $\mB_{red} \subset \N^m$
which can be either empty or satisfies that $\mB_{red} =
\{b_1^{\,\prime},\ldots,b_r^{\,\prime}\}$, where $r \leq n$ and $
b_i^{\,\prime} \in \sum_{j \in \{1,\ldots,r\}\atop j \neq i} \Z\,
b_j^{\,\prime} \setminus \sum_{j \in \{1,\ldots,r\} \atop j \neq i}
\N\, b_j^{\,\prime}$ for all $i \in \{1,\ldots,r\}$. As a
consequence of this construction we have the following result.

\begin{Theorem}\cite[Theorem 2.5]{BGsimplicial} \label{TheoremReduccion}$\IB$ is a complete intersection $\Longleftrightarrow$ either
$\mB_{red} = \emptyset$ or $I_{\mB_{red}}$ is a complete
intersection.
\end{Theorem}

Let $\mA = \{d_1,\ldots,d_n\}$ be a set of $n \geq 2$  positive
integers, in this setting we have that $B_i = {\rm min}\{B \in \Z^+
\, \vert \, B d_i \in \sum_{j \in \{1,\ldots,n\} \atop j \neq i} \Z
d_j\} = \gcd(\mA \setminus \{d_i\})\, /\, \gcd(\mA)$ for all $i \in
\{1,\ldots,n\}$. For $n = 3$ we have the following result, which is
essentially a rewriting of a classical Herzog's result \cite{Herzog}
(see also \cite[Proposition 3]{Watanabe}).

\begin{Proposition}\label{n=3}If $n = 3$, $\IA$ is a complete intersection
$\Longleftrightarrow\ \mA_{red} = \emptyset$.
\end{Proposition}

For $n > 3$ the same characterization does not hold. Nevertheless,
under certain hypothesis, we prove in Proposition \ref{n-1dist} an
analogous characterization for $\IA$ to be a complete intersection.
To present this result we define
\begin{center}$m_i := {\rm min}\left\{b \in \Z^+ \, \vert \, b d_i
\in \sum_{j \in \{1,\ldots,n\} \atop j \neq i} \N d_j\right\}$ for
every $i \in \{1,\ldots,n\}$,
\end{center} and say that a binomial $f \in \IA$ is {\it critical with
respect to $x_i$} if $f = x_i^{m_i} - \prod_{j \in \{1,\ldots,n\}
\atop j \neq i} x_j^{\alpha_j}$. Critical binomials were introduced
by Eliahou in \cite{E} and later studied by Alc\'antar and
Villarreal in \cite{AV}

\begin{Lemma}\label{lemacritico}Let $f_1,\ldots,f_t$ be critical binomials with respect to
$x_{i_1},\ldots,x_{i_t}$ respectively, where $1 \leq i_1 < \cdots <
i_t \leq n$. If $m_i d_i \neq m_j d_j$ for every  $1 \leq i < j \leq
t,$ then there exists a set of binomials $\frak B$ minimally
generating $\IA$ such that $f_1,\ldots,f_t \in \frak B.$
\end{Lemma}
\begin{demo}Let $\{g_1,\ldots,g_{n-1}\}$ be a set of $\mA$-homogeneous binomials
generating $\IA$.  Let $D_1,\ldots, D_{n-1}$ be the $\mA$-degrees of
$g_1,\ldots,g_{n-1}$  respectively and suppose that $f_1 \in \IA$ is
a critical binomial with respect to $x_1$, thus $f_1 = x_1^{m_1} -
\prod_{j \in \{2,\ldots,n\}} x_j^{\alpha_j}$ for some
$\alpha_2,\ldots,\alpha_n \in \N$. Hence $f_1 =
q_{1}g_1+\ldots+q_{n-1}g_{n-1}$\,, where $q_{k}\in \kx$ is an
$\mA$-homogeneous polynomial of degree $m_1 a_1 - D_k\geq 0$ when
$q_{k}\neq 0$ for all $k\in \{1,\ldots,n-1\}$.
 In particular, there exists $k \in \{1,\ldots,n-1\}$ such that $q_k \neq 0$ and the image of
$g_{k}$ under the evaluation morphism which sends $x_j$ to $0$ for
all $j \neq 1$ is equal to $x_1^{D_{k}/d_1}$, we assume that $k =
1$. Thus \,$g_{1}=x_1^{D_1/d_1}-\x^{\,\beta}$\,, where
$\x^{\,\beta}$ is a monomial of $\mA$-degree $D_1$ which does not
involve the variable $x_1$, and hence $D_1 \in \Z^+ d_1\cap
\sum_{j\in\{2,\ldots,n\}} \N\,d_j$. By the definition of $m_1$ we
get the equality $m_1 a_1 =D_1$, which implies that $q_1 \in k$ and
$\{f_1,g_2,\ldots,g_{n-1}\}$ is a minimal set of generators of
$\IA$. Iterating this argument, we get the result. \QED
\end{demo}

\begin{Proposition}\label{n-1dist} If $n-1$ integers from  $m_1 d_1,\ldots,m_n
d_n$ are different, then $\IA$ is a complete intersection
$\Longleftrightarrow$ $\mA_{red} = \emptyset$.
\end{Proposition}
\begin{demo}$(\Leftarrow)$ Follows from Theorem
\ref{TheoremReduccion}. $(\Rightarrow)$ Assume that $m_1 d_1,\ldots,
m_{n-1} d_{n-1}$ are all different, by Lemma \ref{lemacritico} we
get that $\IA = (f_{1},\ldots,f_{n-1})$ where $f_i$ is a critical
binomial with respect to $x_i$ for every $1 \leq i \leq n-1$. We
claim that there exists $j \in \{1,\ldots,n-1\}$ such that $x_j$
does not appear in $f_k$ for all $k \in \{1,\ldots,n-1\} \setminus
\{j\}$. Suppose this claim is false, then we consider the simple
directed graph with vertex set $\{1,\ldots,n-1\}$ and arc set
$\{(j,k)\, \vert\, 1 \leq j,\,k \leq n-1,\ j \neq k$ and $x_j$
appears in $f_k\}$; since the out-degree of every vertex is greater
or equal to one, there is a cycle in the graph. Suppose
 that the cycle is $(1,2,\ldots,k,1)$ with $k \leq n-1$, this means that
$(f_{1},\ldots,f_{k}) \subset (x_{1},\ldots,x_{k})$, so $\IA
\subsetneq H := (x_{1},\ldots,x_{k},f_{k+1},\ldots,f_{n-1})$ but
this is not possible because $\IA$ is prime and $n-1 = {\rm ht}(\IA)
< {\rm ht}(H) \leq n-1$.

Thus there exists $i \in \{1,\ldots,n-1\}$ such that $x_i$ only
appears in $f_i$, suppose $i = 1$. Now we write $\gamma_i := m_i e_i
- \sum_{j \in \{1,\ldots,n\} \atop j \neq i} \alpha_{i,j} e_j \in
\Z^n$ for all $i \in \{1,\ldots,n-1\}$. By \cite[Proposition
2.3]{ElVi}, $\{\gamma_{1},\ldots,\gamma_{n-1}\}$ is a $\Z$-basis for
the kernel of the homomorphism $\tau: \Z^n \longrightarrow \Z$
induced by $\tau(e_j) = d_j$. By definition, $B_1 d_1 = \sum_{j \in
\{2,\ldots,n\}} \beta_j d_j$ for some $\beta_j \in \Z$, so take
$\delta := B_1 e_1 - \sum_{j \in \{2,\ldots,n\}} \beta_j e_j \in
{\rm ker}(\tau)$. Consequently $m_1$ divides $B_1$ and by definition
$B_1$ divides $m_1$, so $B_1 d_1 = m_1 d_1 \in \sum_{j \in
\{2,\ldots,n\}} \N d_j$. Now we have that $\mA_{red} = (\mA
\setminus \{d_1\})_{red}$, and by Lemma \ref{facilred} and
Proposition \ref{red} it follows that $I_{\mA \setminus \{d_1\}}$ is
a complete intersection minimally generated by
$\{f_{2},\ldots,f_{n-1}\}$; repeating the same argument we conclude
that $\mA_{red} = \emptyset$.  \end{demo} \QED

\bigskip Concerning the case of projective monomial curves,
we denote $d := {\rm max}(\mA)$ and $\omA = \{a_1,\ldots,a_{n+1}\}$
where $a_i = (d_i, d - d_i)$ for every $i \in \{1,\ldots,n\}$ and
$a_{n+1} = (0,d)$. Since toric ideals associated to projective
monomial curves is a subfamily of homogeneous simplicial toric
ideals, the following result, which is a particular case of
\cite[Corollary 3.4]{BGsimplicial}, holds.

\begin{Theorem}\label{coralg} $\oIA$ is a complete intersection $\Longleftrightarrow$  $\omA_{red} = \emptyset$.
\end{Theorem}

Moreover, the computation of $\omA_{red}$ is simpler than in the
general case if we take into account the following properties which
are easy to prove.

\begin{Remark}\label{notaproyred}

\begin{itemize}
\item[{\rm (1)}] $B_i = \gcd (\mA \setminus
\{d_i\})\, /\, \gcd(\mA)$ for all $i \in \{1,\ldots,n\}$.
\item[{\rm (2)}] $\forall i \in \{1,\ldots,n\}:  B_i a_i \in \sum_{j \in \{1,\ldots,n+1\} \atop j \neq i} \N  a_j
\Longleftrightarrow B_i d_i = \sum_{j \in \{1,\ldots,n\}\atop j \neq
i} \alpha_j d_j \in \sum_{j \in \{1,\ldots,n\} \atop j \neq i} \N
d_j,$ where $\sum_{j \in \{1,\ldots,n\} \atop j \neq i} \alpha_j
\leq B_i.$ Thus, checking whether $B_i a_i \in \sum_{j \in
\{1,\ldots,n\} \atop j \neq i} \N a_j$ is reduced to determining if
an integer belongs to a subsemigroup of $\N$ with an extra
condition.
\item[{\rm (3)}] If $i = n+1$ or $d_i = d = {\rm
max}(\mA)$, then $B_i a_i \notin \sum_{j \in \{1,\ldots,n\} \atop j
\neq i} \N a_j$.
\end{itemize}
\end{Remark}

In Table \ref{algoritmocurvamonomialproyectiva} we show Algorithm
CI-projective-monomial-curve, an algorithm which receives as input a
set $\mA = \{d_1,\ldots,d_n\} \subset \Z^+$ and determines whether
$\oIA$ is a complete intersection. This algorithm is essentially
obtained in \cite[Theorem 3.6]{CN}.

\begin{table}[!htb]
\centering
\begin{tabular}{|p{11.5cm}|}
\hline

\begin{center}
{\bf Algorithm CI-projective-monomial-curve}
\end{center}
\vspace{-.3cm} $$\begin{array}{ll}
\ Entrada:  \mA = \{d_1,\ldots,d_n\} \subset \Z^+ \\
\ Salida:\mbox{\ {\sc True} o {\sc False}}
\end{array}$$

\begin{center}
\begin{algorithmic}

\STATE $d := {\rm max}\{d_1,\ldots,d_n\}$

\REPEAT

\STATE $\mB := \mA$

\FORALL {$d_i \in \mA \setminus \{d\}$}

\vspace{.2cm}

 \STATE $B_i := \gcd(\mA \setminus
\{d_i\})\, / \, \gcd(\mA)$

\vspace{.2cm}

\IF {$B_i\, d_i = \sum_{d_j \in \mA \, \atop j \neq i} \alpha_j\,
d_j \in \sum _{d_j \in \mA \, \atop j \neq i} \N\, d_j$ and $\sum
\alpha_j \leq B_i$}

\vspace{.2cm}

 \STATE $\mA :=
\mA \setminus \{d_i\}$

\ENDIF

\ENDFOR

\UNTIL $(\mA = \{d\})$ OR $(\mA  = \mB)$

\IF {$\mA = \{d\}$}

\RETURN {\sc True}

\ENDIF

\RETURN {\sc False}

\end{algorithmic}
\end{center}
\\
\hline
\end{tabular}
\caption{Algorithm {\bf CI-projective-monomial-curve}}
\label{algoritmocurvamonomialproyectiva}
\end{table}

\section{Complete intersections and generalized arithmetic sequences}

\smallskip In this section we deal with the cases in which
either $\mA$ is a generalized arithmetic sequence, or $\mA \setminus
\{d_n\}$ is a generalized arithmetic sequence and $d_n$ is a
positive integer. We denote by $\mS$ the semigroup spanned by $\mA$
and assume that {\bf\mathversion{bold}$\mS$ is a numerical
semigroup}, i.e., $\gcd(\mA) = 1$, and  that
{\bf\mathversion{bold}$\mS$ is minimally generated by $\mA$}. Note
that if $\mA$ is a generalized arithmetic sequence and  $\gcd(\mA) =
1$, one can easily check that $\mA$ is a minimal set of generators
of $\mS$ if and only if $n \leq a_1$.

\smallskip As we mentioned in the introduction, Maloo and Sengupta
proved in \cite{MS} that whenever $\mA \setminus \{d_n\}$ is an
arithmetic sequence and  $n \geq 5$, then $\IA$ is not a complete
intersection. To prove this they used the description of a minimal
set of generators of $\IA$ when $\mA \setminus \{d_n\}$ is an
arithmetic sequence obtained by Patil and Singh \cite{PatilSingh}
(see also \cite{Patil} for a shorter proof of the same result). Here
we present a generalization of Maloo and Sengupta's result which
does not require to obtain a description of a minimal set of
generators of $\IA$. More precisely, we prove that whenever $\mA
\setminus \{d_n\}$ is a generalized arithmetic sequence and $n \geq
5$, then $\IA$ is not a complete intersection. In order to prove
this result, we first introduce two results.

\begin{Lemma}\label{noIC}Let $\mA = \{d_1,\ldots,d_n\}$ be a subset of
$\Z^+$. If there exist $1 \leq i < j \leq n$ such that:
\begin{itemize}
\item $m_i d_i \neq m_j d_j$ and
\item $m_i d_i = \sum_{k \in \{1,\ldots,n\} \atop k \neq i} \alpha_k d_k,\, m_j d_j = \sum_{k \in \{1,\ldots,n\} \atop k \neq j} \beta_k d_k$
with $\alpha_k,\beta_k \in \N,\, \alpha_j \neq 0$ and $\beta_i \neq
0$,
\end{itemize} then $\IA$ is not a complete intersection.
\end{Lemma}
\begin{demo}Assume that $\IA$ is a complete intersection and set
$f_i := x_i^{m_i} - \prod_{k \neq i} x_k^{\alpha_k}$ and $f_j :=
x_j^{m_j} - \prod_{k \neq j} x_k^{\beta_k}$. By Lemma
\ref{lemacritico} there exist some binomials $g_3,\ldots,g_{n-1} \in
\IA$ such that $\IA = (f_i, f_j, g_3,\ldots,g_{n-1})$. Therefore
$\IA \subsetneq J := (x_i,x_j,g_3,\ldots,g_{n-1})$, but this is not
possible because $\IA$ is a prime ideal and  $n-1 = {\rm ht}(\IA) <
{\rm ht}(J) \leq n-1$. \QED
\end{demo}

\begin{Proposition}\label{casiaritm}If $n \geq 4$ and $\mA$ contains a generalized arithmetic
sequence with $4$ elements, then $\IA$ is not a complete
intersection.
\end{Proposition}
\begin{demo}Suppose that $\{d_1, d_2, d_3, d_4\} \subset \mA$ is a generalized arithmetic
sequence, i.e., there exists $h \in \Z^+$ such that $\{h d_1, d_2,
d_3, d_4\}$ is an arithmetic sequence. Since $\mA$ is a minimal set
of generators of $\mS$, we have that $m_i > 1$ for all $1 \leq i
\leq n$. Moreover, the equalities $2 d_2 = h d_1 + d_3$ and $2 d_3 =
d_2 + d_4$ prove that $m_2 = m_3 = 2$, hence $\IA$ is not a complete
intersection by Lemma \ref{noIC}. \QED
\end{demo}

\smallskip It is worth pointing out that from \cite[Theorem 2.5]{HKV} one can deduce
a weaker version of Proposition \ref{casiaritm} which states that if
$\mA$ contains an arithmetic sequence with $5$ elements, then $\IA$
is not a complete intersection.

\smallskip From Proposition \ref{casiaritm} one directly derives the following
two corollaries.

\begin{Corollary}\label{casiaritgen5}If $\mA \setminus \{d_n\}$ is a generalized arithmetic sequence and
$n \geq 5$, then $\IA$ is not a complete intersection.
\end{Corollary}

\begin{Corollary}\label{aritgenmenor4}If $\mA$ is a generalized arithmetic sequence and
$n \geq 4$, then $\IA$ is not a complete intersection.
\end{Corollary}

Now we can proceed with the characterizations:

\begin{Theorem} \label{gen}Let $\mA$ be a generalized arithmetic sequence with $n \geq 3$.
Then, $\IA$ is a complete intersection  $\Longleftrightarrow \, n =
3$ and $d_1$ is even.
\end{Theorem}
\begin{demo}By Corollary \ref{aritgenmenor4} it only remains to study when $n = 3$.
Let $h \in \Z^+$ be such that $\{h a_1, a_2, a_3\}$ is an arithmetic
sequence. Since $\gcd(\mA)= 1$, denoting $d := d_3 - d_2$ we have
that $d_2 = h d_1 + d$, $d_3 = h d_1 + 2 d$ and $\gcd\{d_1,d\} = 1$.
We separate two cases, if $d_1$ is even, then $B_2 d_2 = 2 d_2 = h
d_1 + d_3 \in \N \{d_1, d_3\}$ and $\mA_{red} = \emptyset$, thus by
Proposition
 \ref{n=3}, $\IA$ is a complete intersection. If $d_1$ is odd, then
\begin{itemize}
\item $B_1 d_1 = \gcd\{h, d\}\, d_1 < d_2 < d_3$, thus  $B_1 d_1 \not\in \N \{d_2, d_3\},$
\item $B_2 d_2 = d_2 \not\in \N \{d_1, d_3\}$ and
\item $B_3 d_3 = d_3 \not\in \N \{d_1, d_2\}$.
\end{itemize}
So, $\mA_{red} = \{B_1 d_1, d_2, d_3\}$ and again by Proposition
\ref{n=3} we conclude that $\IA$ is not a complete intersection.
\QED
\end{demo}

\medskip

\begin{Remark} Furthermore, whenever $\IA$ is a complete intersection, i.e.,
when $\{d_1,d_2,d_3\}$ is a generalized arithmetic sequence and
$d_1$ is even, we get the following additional information $($see
{\rm Lemma \ref{facilred}} and {\rm Proposition \ref{red})}:

\begin{itemize}
\item $\IA = \left( x_2^2\, -\, x_1^h\, x_3,\,
x_1^{ d_3 / 2}\, -\, x_3^{d_1 / 2} \right)$, with $h \in \Z^+$ such
that $\{h d_1, d_2, d_3\}$ is an arithmetic sequence, i.e., $h = (2
d_2 - d_3) / d_1$.

\item $\g(\mS) = d_1 d_3 / 2 - d_1 +\, d_2\, -\, d_3.$
\end{itemize}

The general formula for the Frobenius number of the semigroup $\N
\mA$ when $\mA$ is a generalized arithmetic sequence can be found in
{\rm \cite[Theorem 3.3.4]{RA}}.
\end{Remark}

\smallskip As a direct consequence of Theorem \ref{gen} we get the already cited results:

\begin{Corollary} \cite[Theorem 3.5]{MS} \label{aritm}
Let $\mA$ be an arithmetic sequence. Then, $\IA$ is a complete
intersection $\Longleftrightarrow$ $n = 2$ or $n = 3$ and $d_1$ is
even.
\end{Corollary}

\begin{Corollary} \cite[Corollary~9]{GR99}
Let $\mA$ be a set of consecutive integers. Then, $\IA$ is a
complete intersection $\Longleftrightarrow$ $n = 2$ or $n = 3$ and
$d_1$ is even.
\end{Corollary}

Concerning the case where $\mA \setminus \{d_n\}$ is a generalized
arithmetic sequence with $n \geq 4$, we have the following:

\begin{Theorem} \label{almost}
Let $\mA \setminus \{d_n\}$ be a generalized arithmetic sequence
with $n \geq 4$. Then, $\IA$ is a complete intersection
$\Longleftrightarrow$ $n = 4$ and one of the following holds:
\begin{enumerate}
\item[{\rm 1.}] $d_1 / \gcd\{d_1,d_2\}$ is even and $\gcd\{d_1,d_2\}  d_4 \in \N
\{d_1,d_2,d_3\}$, or
\item[{\rm 2.}] $d_1, d_4$ are even and $\mC_{red} = \emptyset$ with $\mC =
\{d_1,d_3,d_4\}$\end{enumerate}
\end{Theorem}
\begin{demo}By Corollary \ref{casiaritgen5} it only remains to study when $n = 4$.
By Proposition \ref{red} we have that $\IA$ is a complete
intersection if and only if $I_{\mA'}$ so is, where $\mA' :=
\{d_1,d_2,d_3,B_4 d_4\}$. Moreover, $I_{\mA'} = \IB$ where $\mB :=
\{d_1',d_2',d_3',d_4'\}$ with $d_i' = d_i / B_4$ for $1 \leq i \leq
3$ and $d_4' = d_4$.

We separate two cases, if $d_4' \in \N \{d_1',d_2',d_3'\}$, then by
Lemma \ref{facilred}, $\IA$ is a complete intersection if and only
if $I_{\mB \setminus \{d_4'\}}$ so is. Since $\mB \setminus \{d_4'\}
= \{d_1', d_2', d_3'\}$ is a generalized arithmetic sequence and
$\gcd\{d_1',d_2',d_3'\} = 1$, by Theorem \ref{gen} we get that $\IA$
is a complete intersection if and only if $d_1'$ is even. Assume now
that $d_4' \notin  \N \{d_1',d_2',d_3'\}$ and set $m_i := {\rm
min}\{b \in \Z^+ \, \vert \, b d_i' \in \N (\mB \setminus
\{d_i'\})\}$ for all $1 \leq i \leq 4$, then $m_i \geq 2$ for every
$1 \leq i \leq 4$, in particular $m_2 = 2$. Let us study the
possible values of $m_1$, if $m_1 d_1' = m_2 d_2'$ we set $d_5' :=
\gcd\{d_1',d_2'\} = 1$ and $m_3' := {\rm min}\{b \in \Z^+ \, \vert
\, b d_3' \in \N \{d_4',d_5'\}\} = 1$ and by \cite[Proposition
2.1]{BGS} we have that $\IB$ is not a complete intersection because
$m_3 \neq m_3'$. If $m_1 d_1' = m_3 d_3'$, then necessarily $d_1'$
is even; otherwise $m_1 d_1' = \lcm\{d_1',d_3'\} = d_1' d_3' > d_1'
d_2'$, which contradicts the definition of $m_1$. Then we set $d_5'
:= \gcd\{d_1',d_3'\} = 2$ and $m_i' := {\rm min}\{b \in \Z^+ \,
\vert \, b d_i' \in \sum_{j \in \{2,4,5\} \atop j \neq i} \N d_j\}$
for $i = 2,4$, and observe that $m_2' = 1$ o $m_4' = 1$. Indeed, if
$d_2'$ is even then $m_2' = 1$, if $d_4'$ is even then $m_4' = 1$
and if both are odd, then $m_4' = 1$ if $d_2' < d_4'$ , or $m_2'= 1$
otherwise. Again by \cite[Proposition 2.1]{BGS} we have that $\IB$
is not a complete intersection. If $m_1 d_1' \notin \{m_2 d_2', m_3
d_3'\}$ then $m_1 d_1',\, m_2 d_2'$ and $m_3 d_3'$ are all different
and by Proposition \ref{n-1dist}, $\IB$ is a complete intersection
if and only if $\mB_{red} = \emptyset$. Set $B_i' := \gcd(\mB
\setminus \{d_i'\})$ for $1 \leq i \leq 4$, then $B_3' = B_4' = 1$
and $B_i' d_i' \notin \N (\mB \setminus \{d_i'\})$ for $i = 3,4$.
Let $d', h' \in \Z^+$ be such that $d_2' = h' d_1' + d'$ and $d_3' =
d_2' + d'$. If both $d_1'$ and $d_4'$ are even, then
 $B_2' =  2$ and $2 d_2' = h' d_1'
+ d_3'$, thus $\mB_{red} = \emptyset$ if and only if $\mC_{red} =
\emptyset$, where $\mC = \{d_1,d_3,d_4\}$. In case $d_1'$ or $d_4'$
is odd, we have that $B_2' = 1$ and $B_2' d_2' \notin \N
\{d_1',d_3',d_4'\}$. Concerning $B_1'$, we have that $B_1' =
\gcd\{h',d',d_4'\}$; since $B_1' d_1' \mid h' d_1' < d_2' < d_3'$
then $B_1' d_1' \in \N \{d_2',d_3',d_4'\}$ if and only if $d_4' \mid
B_1' d_1'$. If $d_4' \nmid B_1' d_1'$ then we can conclude that
$\mB_{red} \not= \emptyset$ and $\IA$ is not a complete
intersection. Otherwise, $\mB_{red} = \mB'_{red}$ with $\mB' :=
\{d_2',d_3',d_4'\}$, but $\{d_4', d_2', d_3'\}$ is a generalized
arithmetic sequence and $d_4'$ is odd, then by Proposition \ref{n=3}
and Theorem \ref{gen}, $\mB'_{red} \not= \emptyset$.  \QED
\end{demo}

\medskip

\begin{Remark} Furthermore, whenever $\IA$ is a complete intersection, if we let $h$ be the integer such that
 $\{h d_1, d_2, d_3\}$ is an arithmetic sequence, i.e., $h = (2d_2 - d_3) / d_1$, we have the following
 results $($see {\rm Lemma \ref{facilred}} and {\rm Proposition \ref{red})}:

\begin{enumerate} \item If $d_1 / \gcd\{d_1,d_2\}$ is even,
$\gcd\{d_1,d_2\} d_4  \in \N \{d_1,d_2,d_3\}$ and we take $\beta_1,
\beta_2, \beta_3 \in \N$ such that  $b d_4 = \beta_1 d_1 + \beta_2
d_2 + \beta_3 d_3$, where $b := \gcd\{d_1,d_2\}$, then

\begin{itemize} \item $\IA = \left( x_4^{b} -
x_1^{\beta_1} x_2^{\beta_2} x_3^{\beta_3},\,
 x_2^2 - x_1^h x_3,\, x_1^{d_3 / 2 b}
 - x_3^{d_1 / 2 b} \right)$; and

\item $\g(\mathcal S) = (d_1 d_3\, /\, 2 b) - d_1 + d_2 - d_3 + (b - 1) d_4.$
\end{itemize}

\vskip.2cm \item If both $d_1$ and $d_4$ are even, $\mC_{red} =
\emptyset$ with $\mC := \{d_1,d_3,d_4\}$ and we denote by $\mathcal
S^{\prime}$ the numerical semigroup $\N (d_1 / 2) + \N (d_3 / 2) +
\N (d_4 / 2)$, then

\begin{itemize} \item $\IA = \left(  x_2^2 - x_1^h x_3 \right) + I_{\mathcal C} \cdot
k[x_1,x_2,x_3,x_4]$; and

\item $\g(\mathcal S) = 2\, \g(\mathcal S^{\prime}) + d_2$.
\end{itemize}
\end{enumerate}
\end{Remark}

\section{Complete intersections in Fibonacci sequences.}

Given $p,q \in \Z^+$ with $\gcd\{p,q\} = 1$, our next aim is to
characterize when $\IA$ is a complete intersection where  $\mA =
\{d_1,\ldots,d_n\} \subset \Z^+$ with $d_i = F_{e_i}$ for $1 \leq i
\leq n$, $\{F_n\}$ denotes the $(p,q)$-Fibonacci sequence and
$\{e_1,\ldots,e_n\}$ is a generalized arithmetic sequence. In order
to achieve the characterization we first introduce some basic
properties of linear second order recurrence sequences, focusing
specially on those properties of the $(p,q)$-Fibonacci sequence and
also in the Lucas one. Some of these properties are straight
generalizations of those in \cite{Vaj}, some others can be found in
\cite[Section 5]{HP} and the rest can be easily proved.

\smallskip We denote by $[a]_2$ the $2$-valuation of $a \in \Z^+$, i.e.,
$[a]_2 := {\rm max}\{t \in \N \, : \, 2^t|a\}$, and have the
following result:

\begin{Lemma} \cite[Theorems $\bar{f}$, $\bar{\ell}$ y $\bar{f} \bar{\ell}$ ]{HP} \label{gcdfl}
Let $a,b$ be positive integers and set $d := \gcd\{a,b\}$. Then,
\begin{itemize} \item $\gcd\{F_a,F_b\} = F_d.$

\item $\gcd\{L_a,L_b\} = \left\{ \begin{array}{cl}
L_d &$ if $ [a]_2 = [b]_2 \\
2 &$ if $  [a]_2 \neq [b]_2$, $p,\,q$ are odd and $3|d \\ 2 &$ if
$[a]_2 \neq [b]_2$, $p$ is even and $q$ is odd $\\
1&$ otherwise$
\end{array} \right.$
\item $\gcd\{L_a,F_b\} = \left\{ \begin{array}{cl}
L_d &$ if $[a]_2 < [b]_2 \\
2 &$ if $[a]_2 \geq [b]_2$, $p,\,q$ are odd and $3|d \\ 2 &$ if
$[a]_2 \geq [b]_2$, $p,\, b$
are even and $q$ is odd $\\
1 &$ otherwise $
\end{array} \right.$ \end{itemize}
\end{Lemma}

Observing Lemma \ref{gcdfl} and that every $(p,q)$-Fibonacci
sequence is strictly increasing except for $p = 1$ (because $F_1 =
F_2$), and that every $(p,q)$-Lucas sequence is strictly increasing
except for $p \in \{1,2\}$ (because $L_0 \geq L_1$), we get the
following divisibility properties.

\begin{Corollary}\label{divisible}Let $a,\, b$ be two positive integers, we have the following
properties: \begin{itemize}
\item[{\rm (1)}] If $a \mid b$ then $F_a
\mid F_b$.
\item[{\rm (2)}] If $F_a \mid F_b$ then $a \mid b$, unless if $p = b = 1$ and $a = 2$.
\item[{\rm (3)}] If $a \geq 2$, then $L_a \mid L_b$ if and only if $b / a$ is odd.
\item[{\rm (4)}] If $b$ is even, then $\gcd\{L_a,
L_{a+b}\} = \gcd\{L_a,F_b\}$.
\end{itemize}
\end{Corollary}

\medskip
Denote by $\{U_n\}_{n \in \N}$ any sequence satisfying that $U_{n+2}
= p\, U_{n+1} + q\, U_n$ for all $n \geq 2$ and $U_0,\, U_1 \in \N$
are not both null. In the following result we provide some
properties of theses sequences that we will use in the sequel, all
of them are easy to prove. Property $(4)$ can be found in
\cite[Propositions 5.1 and 5.2]{HP}, properties $(1)$ and $(3)$ are
straight generalizations of the corresponding results for  $p = q =
1$ one can find in \cite{Vaj} and $(2)$ can be easily proved.

\begin{Lemma}\label{propfib}Let $a,b,c,d,e$ be positive integers.
Then,
\begin{itemize}

\item[{\rm (1)}] $U_{a+b} = F_a\,U_{b+1} + q\, F_{a-1}\,U_b$.

\item[{\rm (2)}] If $b \in \N \{a_1,\ldots,a_k\}$ where
$a_1,\ldots,a_k \in \Z^+$, then $F_b \in \N
\{F_{a_1},\ldots,F_{a_k}\}$.

\item[{\rm (3)}] $F_a U_b - F_c U_d = (-1)^e  q^e  (F_{a-e}
U_{b-e}  -  F_{c-e} U_{d-e})$ if $a + b = c + d$ and $e  \leq {\rm
min}\{a,b,c,d\}$.

\item[{\rm (4)}] $L_a = F_{2a} / F_a = F_{a+1} + q F_{a-1}$.

\item[{\rm (5)}] $U_{a+2b} + (-1)^b q^{b}  U_a = L_b  U_{a+b}$.
\end{itemize}
\end{Lemma}

With these basic properties, the following inequalities are easy to
prove.

\begin{Corollary}\label{propfib2}Let $a,b,c,d$ be positive integers. The following
inequalities hold:

\begin{itemize}
\item[{\rm (1)}] $q^b U_a \leq U_{a+2b}$ and equality holds if and only if $a = U_1 = 0$, $b = 1$.

\item[{\rm (2)}] $L_a < F_{a+2}$.

\item[{\rm (3)}] $U_{a+b-2} < F_a U_b < U_{a+b-1} \ $ if $\ a, b \geq 2.$

\item[{\rm (4)}] $F_a U_b  < F_c U_d \ $ if $\ a + b < c + d$.

\item[{\rm (5)}] If $a < c$, $a < d$ and $a + b = c + d$, then $F_a U_b
< F_c U_d\ $ if and only if $a$ is even.

\item[{\rm (6)}] $L_{a+b-1} < L_a L_b < {\rm min}\{L_{a+b+1}, 2\,L_{a+b}\}$.

\item[{\rm (7)}] If $a \leq b$, then $L_a  L_b < L_{a+b}$ if $a$ is odd and
$L_a L_b > L_{a+b}$ if $a$ is even.
\end{itemize}
\end{Corollary}

Let $d_i$ denote the $e_i$-th term of the $(p,q)$-Fibonacci
sequence, where $\{e_1,\ldots,e_n\}$ is a generalized arithmetic
sequence. This is, there exist $h, a, d \in \Z^+$ such that $d_1 :=
F_{a}$ and $d_i := F_{ha + (i-1)d}$ for all $i \geq 2$. As we have
mentioned, we aim at characterizing when $\IA$ is a complete
intersection in terms of the values of $p,q,n,h,a,d$. This objective
is achieved with the following result.

\begin{Theorem} \label{fibonacci}Let $p, q \in \Z^+$ be two relatively prime integers and let
 $\{F_n\}_{n \in \N}$ be the $(p,q)$-Fibonacci sequence. Let $\mA =
\{d_1,\ldots,d_n\}$ be the set with $d_1 := F_a$ and $d_i := F_{ha +
(i-1)d}$ for all $i \in \{2,\ldots,n\}$ where $h, a, d \in \Z^+$ and
$n \geq 3$. Then, $\IA$ is a complete intersection
$\Longleftrightarrow$ one the following holds:
\begin{itemize}
\item[{\rm (a)}] $d$ is odd,
\item[{\rm (b)}] $d \geq a$,
\item[{\rm (c)}] $a = 2d$,
\item[{\rm (d)}] $\gcd\{a,d\} = a-d$ and $a$ is odd, or
\item[{\rm (e)}] $n = 3$ and $2d \mid a$.
\end{itemize}
\end{Theorem}

To prove this theorem we use two previous results, namely Lemma
\ref{pertenece} and Proposition \ref{fib}. The first one
characterizes when $d_3 \in \N \{d_1,d_2\}$ and also proves that
whenever $d_3 \in \N \{d_1,d_2\}$, then $d_i \in \N \{d_1,d_2\}$ for
all $i \geq 3$, which particularly implies that $\IA$ is a complete
intersection. The second one characterizes when $\IA$ is a complete
intersection whenever $d_3 \notin \N \{d_1,d_2\}$.

\begin{Lemma} \label{pertenece} $d_3 \in \N \{d_1, d_2\}$
$\Longleftrightarrow$ $d$ is odd or $F_{2d} \geq \lcm\{d_1,\,F_d\}$.
Moreover, if $d_3 \in \N \{d_1, d_2\}$, then $d_i \in \N
\{d_1,d_2\}$ for all $i \geq 3$.
\end{Lemma}
\begin{demo}If $d$ is odd, by $(5)$ in Lemma \ref{propfib} we get that $d_3 = (q^{d}
F_{ha} / F_a) d_1 + L_d d_2 \in \N \{d_1,d_2\}$ and that $d_{i+2} =
L_d d_{i+1} + q^d d_i$ for all $i \geq 2$, thus $d_i \in \N
\{d_1,d_2\}$ for all $i \geq 3$. Suppose that $d$ is even, set $e :=
\gcd\{d_1,d_2\} = F_{\gcd\{a,d\}}$ and consider the numerical
semigroup $\mathcal S := \N \left\{d_1 / e, d_2 / e \right\}$, its
Frobenius number is $\g(\mathcal S) = ((d_1 d_2) / e - d_1 - d_2) /
e$ (see, e.g. \cite[Theorem 2.1.1]{RA}). If $d \geq a$, then $e
\g(\mathcal S) < (d_1 d_2) / e \leq d_1 d_2 \leq F_d\, d_2 < d_i$
for all $i \geq 3$ and we conclude that
 $d_i \in \N \{d_1,d_2\}$. Note that in this case $F_{2d}
>  d_1 F_d \geq \lcm \{d_1,F_d\}$. Suppose now that  $a >
d$ and let us study the existence of solutions $(x,y) \in \N^2$ to
the linear diophantine equation
\begin{equation}\label{eq} d_1  x + d_2 y = d_3.
\end{equation} By Lemma \ref{propfib} we know that $\left(-q^{d}\, \frac{F_{ha}}{d_1} ,\,
L_d \right) \in \Z^2$ is an integer solution to the equation; thus
the set of integral solutions is
$$\left\{ \left( -q^d\, \frac{F_{ha}}{d_1} + \lambda
\frac{d_2}{e},\, L_d - \lambda \frac{d_1}{e} \right) ,\, \lambda \in
\Z \right\}.$$ We claim that $-q^d \frac{F_{ha}}{d_1}+ \frac{d_2}{e}
> 0$; indeed since $a > d$, then $e = F_{\gcd\{a,d\}} \leq F_{a-d}$
and   $q^d\, e\, F_{ha} = q^{d / 2}\, e\, q^{d / 2}\, F_{ha} < q^{d
/ 2}\, F_{a-d}\, d_2  < d_1\, d_2$. Therefore (\ref{eq}) has a
nonnegative integer solution if and only if $d_1 / e \leq L_d =
F_{2d} / F_d$, which is equivalent to $F_{2d} \geq \lcm\{d_1,
F_d\}$. Finally we have that
$$\g(\mathcal T)\,e \leq \frac{d_1\, d_2}{e} \leq L_d \, d_2 <
F_{d+2}\,d_2 < F_{ha+2d+1} < d_i {\rm \ for\  all\ } i \geq 4$$ and
 $d_i \in \N \{d_1,d_2\}$ for all $i \geq 4$. \QED
\end{demo}

\vskip.5cm Now we study when $\IA$ is a complete intersection
provided $d_3 \notin \N \{d_1,d_2\}$.

\begin{Lemma} \label{fibn=4}
If $2d \mid a$ and $d_3 \not\in \N \{d_1,d_2\}$, then $d_4 \not\in
\N \{d_1,d_2,d_3\}$.
\end{Lemma}
\begin{demo}Assume that there exist $\alpha_1, \alpha_2, \alpha_3
\in \N$ such that $d_4 = \alpha_1  d_1 + \alpha_2  d_2 + \alpha_3
d_3$. Since $d_4 = -q^d  d_2 + L_d  d_3$, we have that $\beta d_3 =
\alpha_1  d_1 + (\alpha_2 + q^d) d_2$ where $\beta := L_d -
\alpha_3$ and the equation $\beta d_3 = x  d_1 + y d_2$ has a
solution $(x,y) \in \N \times \N$. Moreover, the equality $\beta d_3
= - \beta q^d (F_{ha} / d_1) d_1 + \beta L_d  d_2$ yields that the
set of integer solutions to this equation is
$$\left\{ \left(-\beta\, q^d \frac{F_{ha}}{d_1} + \lambda \,
\frac{d_2}{F_d},\, \beta\, L_d - \lambda \, \frac{d_1}{F_d} \right),
\, \lambda \in \Z \right\}.$$ Nevertheless, for all $\lambda
> 0$ we have that
$$\beta\, L_d - \lambda \, \frac{d_1}{F_d} \leq (L_d)^2 -
\frac{d_1}{F_d} \leq (L_d)^2 - \frac{F_{4d}}{F_d} = (L_d)^2 - L_{2d}
L_d < 0;$$ and we can conclude that there is no solution $(x,y) \in
\N \times \N$, a contradiction. \QED
\end{demo}

\begin{Proposition}\label{fib}If $d_3 \not\in \N \{d_1,d_2\}$, then
$\IA$ is a complete intersection $\Longleftrightarrow$ $n = 3$ and
$2d\,|\,a$.
\end{Proposition}

\begin{demo}By Lemma \ref{pertenece}, $d$ is even
and by Lemma \ref{propfib} we have that $d_3 + q^d F_{ha} = L_d d_2$
and $d_{i+1} + q^d d_{i-1} = L_d d_i$ for $3 \leq i \leq n-1$, which
implies that $m_i \leq L_d$ for all $i \in \{2,\ldots,n-1\}$. We
claim that $m_2 d_2,\ldots,m_n d_n$ are all different. Indeed if we
assume that there exist $i,j: 2 \leq i < j \leq n$ such that $m_i
d_i = m_j d_j$, then $m_i d_i \leq L_d d_i = d_{i+1} + q^d d_{i-1} <
2\, d_{i+1}$; which implies that $j = i+1$ and $m_i d_i = d_{i+1}$
and $d_i \mid d_{i+1}$, but this is not possible because $ha +
(i-1)d \nmid ha + id$. Hence by Proposition \ref{n-1dist}, $\IA$ is
a complete intersection if and only if $\mA_{red} = \emptyset$.

Let $\mA'$ be the minimal set of generators of $\N \mA$, then there
exists $k \in \N$ and $4 \leq i_1 < \cdots < i_k \leq n$ such that
$\mA' = \{d_1,d_2,d_3,d_{i_1},\ldots,d_{i_k}\}$. We set
$$B_j' := \frac{\gcd(\mA' \setminus \{d_j\})}{ \gcd(\mA')} {\rm \
for \ all \ } j \in \{1,2,3,i_1,\ldots,i_k\}$$ and have that $B_j' =
1$ for all $j \in \{3,i_1,\ldots,i_k\}$ because $\gcd\{d_1,d_2\} =
\gcd(\mA')$ and $B_1' = \frac{F_{\gcd\{ha,d\}}}{\gcd(\mA')}$.
Moreover $B_1' d_1 \not \in \N (\mA' \setminus \{d_1\})$ because
$B_1'  d_1 \leq F_d  d_1 < d_i$ for all $i \geq 2$. If $[a]_2 \leq
[d]_2$ or there exists $j \in \{1,\ldots,k\}$ such that $i_j$ is
even, we also have that $B_2' = 1$, this implies that $\mA_{red}
\not= \emptyset$ and $\IA$ is not a complete intersection. If $[a]_2
> [d]_2$ and $i_j$ is odd for all $j \in \{1,\ldots,k\}$, then
$\gcd\{a,2d\} = 2 \gcd\{a,d\}$ and $B_2 = L_{\gcd\{a,d\}}$. Suppose
first that $\gcd\{a,d\} < d$, then $B_2 d_2 \leq L_{d-1} d_2 < L_{ha
+ 2d - 2} < d_i$ for all $i \in \{3,\ldots,n\}$ and we claim that
$d_1 \nmid B_2\, d_2$; otherwise we take $\alpha_1,\alpha_2 \in \Z$
such that $d_3 = \alpha_1 d_1 + \alpha_2 d_2$, then $B_2 \mid
\alpha_2$ and we get that $d_1 \mid d_3$, a contradiction. Hence
$\mA_{red} \not= \emptyset$ and $\IA$ is not a complete
intersection. Finally assume that $[a]_2
> [d]_2$ and that $\gcd\{a,d\} = d$ or, equivalently, that  $2 d \mid a$;
if $n \geq 4$, by Lemma  \ref{fibn=4} we have that $d_4 \not\in \N
\{d_1,d_2,d_3\}$, then $i_1 = 4$ and we are in the previous case; if
$n = 3$, then $L_d\, d_2 = q^d\, F_{ha} + d_3 \in \N d_1 + \N d_3$,
$\mA_{red} = \emptyset$ and by Proposition \ref{n=3} we conclude
that $\IA$ is a complete intersection. \QED
\end{demo}

\bigskip

\noindent {\it Proof of Theorem \ref{fibonacci}.} As we proved in
Lemma \ref{pertenece}, $d_3 \in \N \{d_1,d_2\}$ if and only if $d$
is odd or $F_{2d} \geq \lcm\{d_1, F_d\}$. Moreover, $F_{2d} \geq
\lcm\{d_1, F_d\}$ $\Longleftrightarrow$ $F_{2d} F_{\gcd\{a,d\}} \geq
d_1 F_d$, and by (4) and (5) in Corollary \ref{propfib2} this is
equivalent to $\gcd\{a,d\} > a - d$ or $a = 2d$ or $\gcd\{a,d\} = a
- d$ and $a - d$ is odd. Furthermore, $\gcd\{a,d\}
> a - d$ if and only if $d \geq a$. So, we have that $d_3 \in
\N \{d_1,d_2\}$ if and only if $d$ is odd, $d \geq a$, $a = 2d$ or
$\gcd\{a,d\} = a - d$ and $a$ is odd. In this situation we also have
that $d_i \in \N \{d_1,d_2\}$ for all $i \in \{3,\ldots,n\}$ and by
Lemma \ref{facilred} we conclude that $\IA$ is a complete
intersection. If $d_3 \not\in \N \{d_1,d_2\}$, the result follows
from Proposition \ref{fib}.\QED

\bigskip It is worth to mention that the characterization obtained in Theorem \ref{fibonacci}
does not depend on the values of $p,\,q$ or $h$, but only on those
of  $a,\,d$ y $n$.

\medskip
\begin{Remark} Denoting $e := F_{\gcd\{a,d\}}$ and $\mathcal S :=
\sum_{i = 1}^n \N \, (d_i / e)$, whenever $\IA$ is a complete
intersection, we get the following additional information $($see
{\rm Lemma \ref{facilred}} and {\rm Proposition \ref{red})}:

\begin{enumerate} \item If $d$ is odd. Then,
\begin{itemize} \item $\IA = \left(  x_1^{d_2 / e} - x_2^{d_1 / e},
\, x_3 - x_1^{q^d F_{ha} / d_1} x_2^{L_d}, x_4 - x_2^{q^d}
x_{3}^{L_d},\ldots,\, x_n - x_{n-2}^{q^d} x_{n-1}^{L_d} \right)$; y

\item $g\left( \mathcal S  \right) =
\frac{1}{e} \,  \left( d_1 d_2 / e - d_1 - d_2 \right)$.
\end{itemize}
\vskip.3cm
\item If $d \geq a$, or $a = 2d$, or $\gcd\{a,d\} = a - d$ and $a$ is odd. Then,
\begin{itemize} \item $\IA= \left( x_1^{d_2 / e} - x_2^{d_1 / e},
x_3 - x_1^{b_{3,1}} x_2^{b_{3,2}},\, x_4 - x_1^{b_{4,1}}
x_2^{b_{4,2}}, \ldots,x_n - x_1^{b_{n,1}} x_{2}^{b_{n,2}}\right)$,
\\ where
$b_{3,1},\ldots,b_{n,2} \in \Z^+$ satisfy that $b_{i,1}\, d_1 +
b_{i,2}\, d_2 = d_i$ for all $i \in \{3,\ldots,n\}$; and

\item $g\left( \mathcal S \right) =
\frac{1}{e} \,  \left( d_1 d_2 / e - d_1 - d_2 \right)$.
\end{itemize}
\vskip.3cm
\item If $n=3$\,, $2d \mid a$\,, $a \neq 2d$
and $d$ is even. Then $e = F_d$ and
\begin{itemize}

\item $\IA = \left(  x_1^{d_3 / F_{2d}} - x_3^{d_1 / F_{2d}},
\, x_2^{L_d} - x_1^{q^d F_{ha} / d_1} x_3 \right)$; and

\item $g \left( \mathcal S \right) =  \frac{1}{F_d} \left( d_1\, d_3 / F_{2d}
 - d_1  + \left( L_d - 1 \right) d_2 - d_3 \right)$.
\end{itemize}
\end{enumerate}
\end{Remark}

\section{Complete intersections in Lucas sequences.}
Let $p,q \in \Z^+$ be two relatively prime intgers, in this section
$d_i$ denotes the $e_i$-th term of the $(p,q)$-Lucas sequence, where
$\{e_i\}_{1 \leq i \leq n}$ is an arithmetic sequence. This is,
there exist $a, d \in \Z^+$ such that $d_i := L_{a + (i-1)d}$ for
all $i \geq 1$ and we aim at characterizing when $\IA$ is a complete
intersection in terms of the values of $p, q, n, a$ and $d$. This
objective is achieved in Theorem \ref{lucas}. Recall that $[a]_2$
(respect. $[d]_2$) denotes the $2$-valuation of $a$ (respect. $d$).

\begin{Theorem} \label{lucas}Let $p,\,q$ be two relatively prime positive integers and let
$\{L_n\}_{n \in \N}$ be the $(p,q)$-Lucas sequence. Let $\mA =
\{d_1,\ldots,d_n\}$ be the set with $d_i := L_{a + (i-1)d}$ for all
$i \in \{1,\ldots,n\}$ where $a,\,d\in \Z^+$ and $n \geq 3$. Then,
$\IA$ is a complete intersection $\Longleftrightarrow$ one of the
following holds:
\begin{itemize}
\item[{\rm (a)}] $d$ is odd,
\item[{\rm (b)}] $d \geq a$,
\item[{\rm (c)}] $\gcd\{a,d\} = a-d$ and $[d]_2 > [a]_2 \geq 1$,
\item[{\rm (d)}] $n = 3$, $p$ and $a/d$ are odd and $q$ is even,
o
\item[{\rm (e)}] $n = 3$, $3 \nmid d$ and $p,\,q$ and $a / d$ are odd.
\end{itemize}
\end{Theorem}

To prove this theorem we have followed a completely analogous scheme
to the one we used for Theorem \ref{fibonacci} but taking into
account the $(p,q)$-Lucas sequence properties shown in Lemmas
\ref{gcdfl} and \ref{propfib} and Corollaries \ref{divisible} and
\ref{propfib2}. For a sake of brevity we are not including here the
proof of Theorem \ref{lucas}, nevertheless we will state Lemma
\ref{pertenecelucas} and Proposition \ref{luc}, that are,
respectively, the Lucas versions of Lemma \ref{pertenece} and
Proposition \ref{fib} which we have used to prove Theorem
\ref{lucas}. Lemma  \ref{pertenecelucas} provides a characterization
of when $d_3 \in \N \{d_1,d_2\}$ and states that whenever $d_3 \in
\N \{d_1,d_2\}$, then $d_i \in \N \{d_1,d_2\}$ for all $i \geq 3$,
which, in particular, implies that $\IA$ is a complete intersection.
Proposition \ref{luc} characterizes when $\IA$ is a complete
intersection whenever $d_3 \notin \N \{d_1,d_2\}$.

\begin{Lemma}\label{pertenecelucas} $d_3 \in \N \{d_1,
d_2\}$ if and only if $d$ is odd or $F_{2d} \geq \lcm\{d_1,F_d\}$.
Moreover, whenever $d_3 \in \N \{d_1,d_2\}$, then $d_i \in \N \{d_1,
d_2\}$ for all $i \geq 3$.
\end{Lemma}

\vskip.5cm Now we study when $\IA$ is not a complete intersection
provided $d_3 \notin \N \{d_1,d_2\}$.

\begin{Lemma} \label{lucasn=4}
If $a / d$ is odd, $\gcd\{d_1, d_2\} = 1$ and $d_3 \not\in \N \{d_1,
d_2\}$, then $d_4 \not\in \N \{d_1, d_2, d_3\}$.
\end{Lemma}

\begin{Proposition}\label{luc}If $d_3 \not\in \N \{d_1,d_2\}$,
then $\IA$ is a complete intersection $\Longleftrightarrow$  $n =
3$, $a/d$ is odd and $\gcd\{d_1,d_2\} = 1$.
\end{Proposition}

\medskip
\begin{Remark} Denoting $e := \gcd(\mA)$ and $\mathcal S :=
\sum_{i = 1}^n \N \, (d_i / e)$, whenever $\IA$ is a complete
intersection, we get the following additional information $($see
{\rm Lemma \ref{facilred}} and {\rm Proposition \ref{red})}:

\begin{enumerate}
\item If $d$ is odd. Then,
\begin{itemize}
\item $\IA = \left(  x_1^{d_2 / e} - x_2^{d_1 / e},
\, x_3 - x_1^{q^d} x_2^{L_d}, \ldots,\, x_n - x_{n-2}^{q^d}\,
x_{n-1}^{L_d} \right)$; and

\item $g\left( \mathcal S\right) =
\frac{1}{e} \,  \left( d_1 d_2 / e - d_1 - d_2 \right)$.
\end{itemize}
\vskip.3cm
\item If $d \geq a$, or $\gcd\{a,d\} = a-d$ and $[d]_2 > [a]_2 \geq 1$. Then,
\begin{itemize}

\item $\IA =  \left(x_1^{d_2 / e} - x_2^{d_1 / e},
x_3 - x_1^{b_{3,1}} x_2^{b_{3,2}},\, x_4 - x_1^{b_{4,1}}
x_2^{b_{4,2}}, \ldots,x_n - x_1^{b_{n,1}} x_{2}^{b_{n,2}}\right)$,
with $b_{i,1}, b_{i,2} \in \N$ such that $b_{i,1} d_1 + b_{i,2} d_2
= d_i$ for all $i \in \{3,\ldots,n\}$; and

\item
$g\left( \mathcal S \right) = \frac{1}{e} \, \left( d_1 d_2 / e -
d_1 - d_2 \right)$.
\end{itemize}
\vskip.3cm
\item If $n=3$, $d$ is even, $a / d$ and $p$ are odd and, either $q$ is even, or $3
\nmid d$. Then $e = 1$ and

\begin{itemize}

\item $\IA = \left(  x_1^{d_3 / L_d} - x_3^{d_1 / L_d},
\, x_2^{L_d} - x_1^{q^d}\, x_3 \right)$; and

\item $\g(\mathcal S) = d_1 d_3 / L_d - d_1 + (L_d -1)d_2 - d_3$.
\end{itemize}
\end{enumerate}
\end{Remark}

\medskip
Theorem \ref{lucas} depends on the values of $a, d, n$ and on the
parity of $p$ and $q$, in contrast to the corresponding result for
the Fibonacci sequence where the values of $p,q$ and even of $h$ do
not play any role in the result (see Theorem \ref{fibonacci}). This
dependence of "all" initial values gives the insight that a more
general result where the complete intersection property for $\IA$ is
characterized, where $d_i = L_{e_i}$ with $\{e_1,\ldots,e_n\}$ a
generalized arithmetic sequence, the value of $h$ such that $\{h
e_1,e_2,\ldots,e_n\}$ is an arithmetic sequence is relevant. The
following example shows the relevance of the value of $h \in \Z^+$,
and gives a taste of the difficulty that might have the more general
case in which $e_1,\ldots,e_n$ is a generalized arithmetic sequence.

\begin{example} Let $\mA_h$ be the set $\{L_{5}, L_{h 5 + d}, L_{h 5 + 2d}\},$
where $\{L_n\}_{n \in \N}$ is the $(1,1)$-Lucas sequence and $h \in
\Z^+$. For $h = 1$ and $h = 3$, the sets $\mA_1 = \{L_5 = 11, L_6 =
18, L_7 = 29\}$ and $\mA_3 = \{L_5 = 11, L_{16} = 2207, L_{17} =
3571\}$ determine complete intersection toric ideals. Nevertheless,
for $h = 2$ and $h = 4$, the sets $\mA_2 = \{L_5 = 11, L_{11} = 199,
L_{12} = 322\}$ and $\mA_4 = \{L_5 = 11 , L_{21} = 24476, L_{22} =
39603\}$ determine two non complete intersection toric ideals.
\end{example}

\section{Complete intersections in certain projective monomial curves}

In this section we denote $\mA = \{d_1,\ldots,d_n\} \subset \Z^+$
and $d := {\rm max} \{d_1,\ldots,d_n\}$ and consider $$\omA :=
\{a_1, \ldots, a_{n-1}, a_n, a_{n+1}\} \subset \N^2,$$ where $a_i :=
(d_i,d - d_i)$ for every $i \in \{1,\ldots,n\}$ and $a_{n+1} :=
(0,d)$.

\smallskip The objective of this section is to characterize when $\oIA$ is a complete intersection
in terms of the set $\mA$, where $\mA$ belongs to any of the
families studied in the previous sections. For this purpose we use
Algorithm CI-projective-monomial-curve of Table
\ref{algoritmocurvamonomialproyectiva}.

\medskip We begin studying when $\oIA$ is a complete intersection, where either $\mA$
or $\mA \setminus \{d_n\}$ is a generalized arithmetic sequence and
we assume without loss of generality that
{\mathversion{bold}$\gcd(\mA) = 1$}.

\begin{Theorem}\label{genproy} Let $\mA = \{d_1,\ldots,d_n\}$ be a
generalized arithmetic sequence with $n \geq 3$. Then, $\oIA$ is a
complete intersection $\Longleftrightarrow$ $n = 3$, $\mA$ is an
arithmetic sequence and $d_1$ is even.
\end{Theorem}
\begin{demo}$(\Rightarrow)$ Set $h \in \Z^+$ such that $\{h d_1, d_2, \ldots, d_n\}$ is
an arithmetic sequence. For every $i \in \{3,\ldots,n-1\}$ we have
that $B_i = 1$ and $B_1 = \gcd\{d_1,h\}$, thus $B_1 d_1 \leq h d_1 <
d_2 < \cdots < d_n$  and $B_1 a_1 \notin \sum_{j=2}^{n+1} \N a_j$.
If $n \geq 4$ or $d_1$ is odd, it follows that $B_2 = 1$ and
$\omA_{red} \not= \emptyset$, thus $\oIA$ is not a complete
intersection by Theorem \ref{coralg}. Suppose now that $n = 3$ and
$d_1$ even, then $B_2 = 2$. If $h \geq 2$, then $B_2 d_2 \not=
\alpha_1 d_1 + \alpha_3 d_3$ with $\alpha_1,\alpha_3 \in \N$ and
$\alpha_1 + \alpha_3 \leq 2$ and, again by Theorem \ref{coralg}, we
have that $\oIA$ is not a complete intersection.

$(\Leftarrow)$ We have that $B_2 = 2$ and $2 a_2 = a_1 + a_3$, then
by Theorem \ref{coralg}, $\oIA$ is a complete intersection. \QED
\end{demo}

\medskip

\begin{Remark} Whenever $\oIA$ is a complete intersection, i.e., when $n = 3$, $\mA$
is an arithmetic sequence and $d_1$ is even, we get the following
minimal set of generators of the toric ideal {\rm (}see {\rm Lemma
\ref{facilred}} and {\rm Proposition \ref{red})}:
$$\oIA = \left( x_2^2\, -\, x_1  x_3,\, x_1^{ d_3 / 2}\, -\,
x_3^{d_1 / 2} x_4^{d_2 - d_1} \right).$$
\end{Remark}

\bigskip

Concerning when  $\mA \setminus \{d_n\}$ is a generalized arithmetic
sequence, we can assume without loss of generality that
$\{d_n,d_2,\ldots,d_{n-1}\}$ is not an arithmetic sequence.
Otherwise we take $d_1' := d_n$, $d_i' = d_i$ for all $2 \leq i \leq
n-1$ and $d_n' = d_1$ and have that  $\mA \setminus \{d_n'\} =
\{d_1',\ldots,d_{n-1}'\}$ is an arithmetic sequence.

\begin{Theorem} \label{almostN=4proy}
Let $\mA = \{d_1,\ldots,d_n\}$ be a set such that $\mA \setminus
\{d_n\}$ is a generalized arithmetic sequence with $n \geq 4$. Then,
$\oIA$ is a complete intersection $\Longleftrightarrow$ $n = 4$,
$\{d_1,d_2,d_3\}$ is an arithmetic sequence and one of the following
holds:
\begin{enumerate}
\item[{\rm 1.}] $d_1\, /\, \gcd\{d_1,d_2\}$ is even
and $\gcd\{d_1,d_2\}\, d_4 = \beta_1 d_1 + \beta_2 d_2 + \beta_3
d_3$ with $\beta_1 + \beta_2 + \beta_3 \leq \gcd\{d_1,d_2\}$, or
\item[{\rm 2.}] $d_1, d_4$ are even and
${\mathcal C}^{\star}_{red} = \emptyset$ with $\mathcal C =
\{d_1,d_3,d_4\}$.
\end{enumerate}
\end{Theorem}
\begin{demo}Let $h \in \Z^+$ be such that $\{hd_1,d_2,\ldots,d_{n-1}\}$ is an arithmetic sequence.
We divide the proof in two parts, if $B_n a_n \in \sum_{j \in
\{1,\ldots,n-1,n+1\}} \N a_j$ or equivalently if $B_n d_n = \sum_{j
= 1}^{n-1} \beta_j d_j$ with $\sum_{j = 1}^{n-1} \beta_j \leq B_n$
(see Remark \ref{notaproyred}), then by Proposition \ref{red} and
Lemma \ref{facilred}  it follows that $\oIA$ is a complete
intersection if and only if lo is $I_{\omA \setminus \{a_n\}}$. It
is easy to check that $I_{\omA \setminus \{a_n\}} = I_{\omA_1}$ with
$\mA_1 = \{d_1/B_n, \ldots, d_{n-1} / B_n\}$. Moreover $\gcd(\mA_1)
= 1$ and $\mA_1$ is a generalized arithmetic sequence, then by
Theorem \ref{genproy} we conclude that $I_{\omA_1}$ is a complete
intersection if and only if $n = 4$, $d_1 / B_4$ is even and $\{d_1,
d_2, d_3\}$ is an arithmetic sequence; note also that $B_4 =
\gcd\{d_1,d_2\}$.

Suppose now that $B_n a_n \notin \sum_{j \in \{1,\ldots,n-1,n+1\}}
\N a_j$. We have that $B_i = 1$ for $i \in \{3,\ldots,n-1\}$ and
then $B_i a_i \notin \sum_{j \in \{1,\ldots,n+1\} \atop j \neq i} \N
a_j$.  Let us study the values of $B_1$ and $B_2$. Set $h,r \in
\Z^+$ such that $d_i = h d_1 + (i-1)r$ for all $i \in
\{2,\ldots,n-1\}$, we have that $B_1 = \gcd\{h, r, d_n\}$, and then
$B_1 d_1 < d_2 < \cdots < d_{n-1}$, so $B_1 d_1 \in \N
\{d_2,\ldots,d_n\}$ if and only if $d_n \mid B_1 d_1$. However, if
$d_n \mid B_1 d_1$, by Proposition \ref{red} and Lemma
\ref{facilred} we have that $\oIA$ is a complete intersection if and
only if $\oIB$ so is, where $\mB := \{d_2, \ldots, d_n\}$, but $d_n
\mid B_1 d_1 \mid h d_1$ and $\mB$ is a generalized arithmetic
sequence; therefore, by Theorem \ref{genproy} this implies that $n =
4$ and $\{d_4, d_2, d_3\}$ is an arithmetic sequence, what we had
assumed not to happen. Concerning $B_2$, if $n \geq 5$, $d_1$ is odd
or $d_4$ is odd then $B_2 = 1$ and $B_2 a_2 \notin \sum_{j \in
\{1,\ldots,n+1\} \atop j \neq 2} \N a_j$. Otherwise, i.e., if $n =
4$ and $d_1$ and $d_4$ are even, we have that $B_2 = 2$ and set $\mC
:= \{d_1,d_3,d_4\}$. Moreover, $2 a_2 \in \N \{a_1,a_3,a_4,a_5\}$ if
and only if $2 d_2 = \alpha_1 d_1 + \alpha_3 d_3 + \alpha_4 d_4$
with $\alpha_1 + \alpha_3 + \alpha_4 \leq 2$, and this can only
happen if:
\begin{itemize}
\item[{\rm (a)}] $2 d_2 = d_1 + d_3$,
\item[{\rm (b)}] $2 d_2 = d_1 + d_4$, or
\item[{\rm (c)}] $2 d_2 = d_3 + d_4$.
\end{itemize}
If (a) holds we have that $\{d_1, d_2, d_3\}$ is an arithmetic
sequence and $\oIA$ is a complete intersection if and only if so is
$I_{\mC^{\star}}$. In (b) we have that $2 d_2 = h d_1 + d_3 = d_1 +
d_4$ and we deduce that $(h-1) d_1 + d_3 = d_4$, thus
 one derives that ${\mathcal C}^{\star}_{red} \not=
\emptyset$ and $\oIA$ is not a complete intersection. Indeed,
denoting $B_i' := {\rm min}\{b \in \Z^+ \, \vert \, b a_i \in
\sum_{j \in \{1,3,4\} \atop j \neq i} \Z a_j\}$ for all $i \in
\{1,3,4\}$, from the previous equality it follows that $B_3' = B_4'
= 1$ and that $B_1' \leq h-1$, therefore $B_1' d_1 \notin \N
\{d_3,d_4\}$, due to $B_1' d_1 < d_3 < d_4$. In (c) we have that
$\{d_4,d_2,d_3\}$ is an arithmetic sequence, which we had assumed to
not to happen.

In the rest of cases we get that $B_i a_i \notin \N (\omA \setminus
\{a_i\})$ for all $i \in \{1,\ldots,n\}$, and then $\omA_{red} \not=
\emptyset$ and $\oIA$ is not a complete intersection.  \QED
\end{demo}

\begin{Remark} Moreover, whenever $\oIA$ is a complete intersection, we get the following
minimal sets of generators of $\oIA$ depending on the cases
 $($see {\rm Lemma \ref{facilred}} and
{\rm Proposition \ref{red}}$)$:

\begin{enumerate} \item Set $b = \gcd\{d_1,d_2\}$, if $n = 4$, $\{d_1,d_2,d_3\}$ is an arithmetic sequence, $d_1 / b$ is even and $b d_4 =
\beta_1 d_1 + \beta_2 d_2 + \beta_3 d_3$ with $\beta_1 + \beta_2 +
\beta_3 \leq b$, then
$$\oIA = \left( x_4^{b} - x_1^{\beta_1} x_2^{\beta_2} x_3^{\beta_3}
x_5^{b - \sum \beta_i},\,
 x_2^2 - x_1 x_3,\, x_1^{d_3 / 2 b}
 - x_3^{d_1 / 2 b} x_5^{(d_2 - d_1) / b} \right).$$

\item If $n = 4$, $\{d_1,d_2,d_3\}$ is an arithmetic sequence, $d_1, d_4$ are even
and ${\mathcal C}^{\star}_{red} = \emptyset$ with $\mathcal C =
\{d_1,d_3,d_4\}$, then
$$\oIA = \left(  x_2^2 - x_1 x_3 \right) + I_{\mathcal C^{\star}} \cdot
k[x_1,x_2,x_3,x_4,x_5].$$
\end{enumerate}
\end{Remark}

\bigskip We finish this section studying when $\oIA$ is a complete intersection in the following families, where $\mA =
\{d_1,\ldots,d_n\}$, $n \geq 3$, $p,q \in \Z^+$ are relatively prime
and either
\begin{itemize}
\item $\mA$ consists of terms the $(p,q)$-Fibonacci sequence whose indices are
a generalized arithmetic sequence, i.e., there exist $h, a, d \in
\Z^+$ such that $d_1 = F_{a}$, $d_i = F_{ha + (i-1)d}$ for all $i
\geq 2$, or

\item $\mA$ consists of terms of the  $(p,q)$-Lucas sequence whose indices are an arithmetic
sequence, i.e., there exist $a, d \in \Z^+$ such that $d_1 = L_{a}$,
$d_i = L_{ha + (i-1)d}$ for all $i \geq 2$.
\end{itemize}

In both cases we characterize when $\oIA$ is a complete intersection
by means of the input data; which are the values of $p, q, n, a$ and
$d$ (and $h$ for the Fibonacci sequence).

\begin{Theorem} \label{fibonacciproy}Let $p, q \in \Z^+$ be relatively prime and
let $\{F_n\}_{n \in \N}$ be the $(p,q)$-Fibonacci sequence. Set $\mA
= \{d_1,\ldots,d_n\}$ with $d_1 := F_a$ y $d_i := F_{ha + (i-1)d}$
for all $i \in \{2,\ldots,n\}$ where $h,\,a,\,d\in \Z^+$ and $n \geq
3$. Then, $\oIA$ is a complete intersection $\Longleftrightarrow$ $n
= 3$, $h = 1$, $d$ is even and $2d \mid a$.
\end{Theorem}
\begin{demo}We begin by observing that $B_i = 1$ for all $i \geq 3$
(see Lemma \ref{gcdfl}) and $B_1 = F_{\gcd\{ha,d\}} /
F_{\gcd\{a,d\}}$ and thus $B_1 d_1 \leq F_d d_1 < d_2 < \cdots <
d_n$, which implies that $B_1 a_1 \notin \sum_{j = 2}^{n+1} \N a_j$.
Let us study the possible values of $B_2$. If $n \geq 4$ or $[d]_2
\geq [a]_2$, then $B_2 = 1$, and if $n = 3$ and $[d]_2 < [a]_2$ then
$B_2 = F_{\gcd\{a,2d\}} / F_{\gcd\{a,d\}} = L_{\gcd\{a,d\}}$. If $d
\nmid a$, then $B_2 d_2 \leq L_{d-1} d_2 < L_{ha+2d-2} < d_3$ and
$B_2 d_2 \in \N d_1 + \N d_3$ if and only if $B_2 d_2 = \alpha d_1$,
but in this case $\alpha
> B_2$ and, by Remark \ref{notaproyred}, $B_2 a_2 \notin \N a_1 + \N a_3 + \N a_4$.
Suppose now that $d \mid a$, if $d$ is odd then $L_d d_2 = - q^d
(F_{ha}/d_1) d_1 + d_3 < d_3$ and again we have that  $B_2 d_2 \in
\N d_1 + \N d_3$ if and only if $B_2 d_2 = \alpha d_1$, and thus
$B_2 a_2 \notin \N a_1 + \N a_3 + \N a_4$. Finally if $d$ is even we
have the inequality $B_2 d_2 = q^d (F_{ha}/d_1) d_1 + d_3 < 2 d_3$,
hence whenever $B_2 d_2 = \alpha_1 d_1 + \alpha_3 d_3$ with
$\alpha_1,\alpha_3 \in \N$, then $\alpha_3 < 2$; and if $\alpha_3 =
0$, then $\alpha_1 > B_2$. As a consequence, by Remark
\ref{notaproyred}, if follows that $B_2 a_2 \in \N a_1 + \N a_3 + \N
a_4$ if and only if $q^d (F_{ha} / d_1) + 1 \leq L_d$. If $h
> 1$, then $F_{ha}/d_1 \geq F_{2a}/ F_a = L_a \geq L_d$ and if $h =
1$ it follows that $L_d = F_{d+1} + q F_{d-1} > q^{d/2}\, F_1 + q\,
q^{d/2 - 1}\, F_1 = q^d$, and the result follows. \qed
\end{demo}

\medskip

\begin{Remark}Whenever $\oIA$ is a complete intersection, i.e., when $\mA = \{d_1,d_2,d_3\}$ with $d_1 = F_a$\,,
$d_2 = F_{a+d}$\,, $d_3 = F_{a+2d}$\,, $d$ is even and $2d \mid a$,
we get the following minimal set of generators of $\oIA$ $($see {\rm
Lemma \ref{facilred}} and {\rm Proposition \ref{red}}$):$

\begin{center}
$\oIA = \left(x_1^{ d_3 / F_{2d}}  -  x_3^{d_1 / F_{2d}}\, x_4^{(d_3
- d_1) / F_{2d}},\,  x_2^{L_d}  -  x_1^{\,q^d}\, x_3\, x_4^{\, L_d -
q^d - 1} \right)$. \end{center}
\end{Remark}

\vspace{.4cm}

\begin{Theorem} \label{lucasproy}Let $p, q \in \Z^+$ be relatively prime and
let $\{L_n\}_{n \in \N}$ be the $(p,q)$-Fibonacci sequence. Set $\mA
= \{d_1,\ldots,d_n\}$ with $d_i := L_{a + (i-1)d}$ for all $i \in
\{1,\ldots,n\}$ where $a,\,d\in \Z^+$ and $n \geq 3$. Then, $\oIA$
is a complete intersection $\Longleftrightarrow$ $n = 3$, $d$ is
even, $p$ and $a/d$ are odd and, either $q$ is even or $3 \nmid d$
\end{Theorem}
\begin{demo}We begin by observing that $B_i = 1$ for all $i \geq 3$
(see Lemma \ref{gcdfl}). Let us study the possible values of
 $B_2$, if $n \geq 4$ or $[d]_2
\neq [a]_2$, then $B_2 = 1$, and if $n = 3$ and $[d]_2 = [a]_2$ then
$B_2 = \gcd\{d_1,d_3\}\, / \, \gcd(\mA) = L_\gcd\{a,d\} \, /\,
\gcd(\mA)$ and $\gcd(\mA) \in \{1,2\}$. If $\gcd(\mA) = 2$ or
$\gcd\{a,d\} < d$, then $B_2 d_2 < d_3$ and thus $B_2 a_2 \notin \N
a_1 + \N a_3 + \N a_4$. Otherwise, if $\gcd(\mA) = 1$ and
$\gcd\{a,d\} = d$, then $B_2 d_2 = L_d d_2 = (-1)^d q^d d_1 + d_3$,
if $d$ is odd then $L_2 d_2 < d_3$ and again we have that  $B_2 a_2
\notin \N a_1 + \N a_3 + \N a_4$; otherwise $B_2 d_2 = q^d d_1 +
d_3$ and $q^d + 1 \leq L_d$. From here we deduce that  $\omA_{red} =
\emptyset$ if and only if $n = 3$, $[d]_2 = [a]_2$, $\gcd\{a,d\}=
d$, $\gcd\{d_1,d_2,d_3\} = 1$ and $d$ is even, by Lemma \ref{gcdfl}
the result follows. \qed
\end{demo}

\medskip

\begin{Remark} Whenever $\oIA$ is a complete intersection, i.e., when $\mA = \{d_1,d_2,d_3\}$ with $d_1 = L_a$\,,
$d_2 = L_{a+d}$\,, $d_3 = L_{a+2d}$\,, $d$ is even, $p$ and $a/d$
are odd and, either $q$ is even or $3 \nmid d$, we get the following
minimal set of generators of $\oIA$ $($see {\rm Lemma
\ref{facilred}} and {\rm Proposition \ref{red}}$):$

\begin{center}
$\oIA = \left( x_1^{d_3 / L_d} - x_3^{d_1 / L_d}\,  x_4^{(d_3 - d_1)
/ L_{d}}, \, x_2^{L_d} - x_1^{q^d}\, x_3 \, x_4^{\, L_d - q^d - 1}
\right)$. \end{center}
\end{Remark}

\bibliographystyle{plain}

\begin{thebibliography}{10}
\bibitem{AV}{A. Alc\'antar and R.~H. Villarreal,
Critical binomials of monomial curves, Comm. Algebra {\bf 22(8)}
(1994) 3037--3052.}

\bibitem{BGRV}
{I. Bermejo, Ph. Gimenez, E. Reyes and R.~H. Villarreal, Complete
intersections in affine monomial curves, Bol. Soc. Mat. Mexicana 3a.
Serie {\bf 11~(2)} (2005) 191--204.}

\bibitem{BGS}
{I. Bermejo, I. Garc\'{\i}a-Marco and J.~J. Salazar-Gonz\'alez, An
algorithm for checking whether the toric ideal of an affine monomial
curve is a complete intersection, J. Symbolic Computation {\bf 42}
(2007) 971--991.}



\bibitem{BGlib2} {I. Bermejo and I. Garc\'{\i}a-Marco, {\tt
cisimplicial.lib}. A distributed {\sc Singular}~3-1-6 library for
determining whether a simplicial toric ideal is a complete
intersection (2012).}


\bibitem{BGsimplicial}
{I. Bermejo and I. Garc\'{\i}a-Marco, Complete intersections in
simplicial toric varieties, J. Symbolic Computation (2014), to appear.}


\bibitem{CN}{M.~P. Cavaliere and G. Niesi, Sulle curve monomiali proiettive
intersezione completa, Bollettino U.M.I. Algebra e Geometria Serie
VI, Vol III-D,\, N.1\, (1984), 189--200.}


\bibitem{DGPS} W. Decker, G.-M. Greuel, G. Pfister and H. Schoenemann, {\sc Singular}~3-1-6,
a Computer Algebra System for Polynomial Computations, Center for
Computer Algebra, University of Kaiserslautern. Available at {\tt
http://www.singular.uni-kl.de} (2012).


\bibitem{D} C. Delorme, Sous-mono\"{\i}des d'intersection compl${\rm
\grave{e}}$te de $\mathbb{N}$, Ann. Sci. \'Ecole Norm. Sup. {\bf 9}
(1976) 145--154.

\bibitem{E} {S. Eliahou, {\it Courbes monomiales et alg\`ebre de
Rees symbolique}, PhD Thesis, Universit\'e de Gen\`eve, 1983.}


\bibitem{ElVi}{S. Eliahou and R.~H. Villarreal, On systems of
binomials in the ideal of a toric variety, Proc. Amer. Math. Soc.
{\bf 130} (2002) 345--351.}

\bibitem{Fel} L. Fel, Symmetric numerical semigroups generated by
 Fibonacci and Lucas triples, Integers {\bf 9}  (2009), A9, 107--116

\bibitem{GR99} P. A. Garc\'{\i}a-S\'{a}nchez and J.~C. Rosales, Numerical semigroups generated by intervals, Pacific
J. Math. {\bf 191~(1)} (1999) 75--83.


\bibitem{Herzog}{J. Herzog, Generators and relations of abelian
semigroups and semigroup rings, Manuscripta Math. \textbf{3} (1970)
175--193.}

\bibitem{HP}{P. Hilton and J. Pedersen, Fibonacci and Lucas Numbers in
Teaching and Research, J. Math. Informatique {\bf 3}, 36--57,
1991-1992.}

\bibitem{HKV}{A. Hosry, Y. Kim and J. Validashti, On the equality of
ordinary and symbolic powers of ideals, J. Commut. Alg. {\bf 4}
(2012), no. 2, 281--292.}

\bibitem{MS}
{A.~K. Maloo and I. Sengupta, Criterion for complete intersection of
certain monomial curves, in: Advances in algebra and geometry
(Hyderabad, 2001), Hindustan Book Agency, New Delhi, 2003, pp.
179--184.}


\bibitem{Patil} {D.~P. Patil, Minimal sets of generators for the relation ideals of certain monomial curves,
Manuscripta Math. {\bf 80}, Number 1 (1993), 239--248.}

\bibitem{PatilSingh} {D.~P. Patil and B. Singh, Generators for the derivation modules and the relation ideals of certain
curves, Manuscripta Math. {\bf 68}, Number 1 (1990), 327--335.}


\bibitem{RA}{J.~L. Ram\'{\i}rez Alfons\'{\i}n, The Diophantine Frobenius Problem,
 Oxford Lecture Series in Mathematics and its Applications {\bf 30}, Oxford University Press, 2005.}


\bibitem{monalg}{R.~H. Villarreal, Monomial Algebras,
Monographs and Textbooks in Pure and Applied Mathematics {\bf 238},
Marcel Dekker, New York, 2001.}

\bibitem{Sturm} {B.~Sturmfels, {\it Gr\"{o}bner Bases and Convex Polytopes\/},
University Lecture Series {\bf 8}, American Mathematical Society,
Providence, RI, 1996.}

\bibitem{Vaj}{S. Vajda, Fibonacci and Lucas numbers, and the Golden Section:
Theory and Applications, New York: Halsted Press, John Wiley \&
Sons, 1989.}

\bibitem{Watanabe}{K. Watanabe, Some examples of
one dimensional Gorenstein domains, Nagoya Math. J. {\bf 49} (1973)
101--109.}

\end{thebibliography}

\end{document}